\documentclass[preprint,12pt]{elsarticle}

\usepackage{amsmath}
\usepackage{amssymb}
\usepackage{amsfonts}
\usepackage{amsthm}
\usepackage{enumitem}
\usepackage{fullpage}
\usepackage{hyperref}

\newtheorem{theorem}{Theorem}[section]
\newtheorem{lemma}[theorem]{Lemma}
\newtheorem{proposition}[theorem]{Proposition}
\newtheorem{corollary}[theorem]{Corollary}
\newtheorem{problem}[theorem]{Problem}

\theoremstyle{definition}
\newtheorem{definition}[theorem]{Definition}
\newtheorem{remark}[theorem]{Remark}
\newtheorem{example}[theorem]{Example}
\journal{~}

\begin{document}

\begin{frontmatter}

\title{Finite groups with no two normal subgroups of the same order}

\author{Angsuman Das}
\ead{angsuman.maths@presiuniv.ac.in}

\address{Department of Mathematics, \\Presidency University, \\86/1 College Street, Kolkata 700073, India\\}

\begin{abstract}
A finite group is said to have property $\mathcal{D}$ if distinct normal subgroups have distinct orders. Cyclic groups and simple groups are obvious examples. We show that $\mathcal{D}$ is implied by property $\mathcal{P}$ (for every prime $p$, cyclic $p$-subgroups of equal order are conjugate), and hence by each of the stronger conjugacy conditions on subgroups of equal order considered in the literature; the symmetric group $S_4$ shows that the converse fails. Our main result is that a supersolvable group has property $\mathcal{D}$ if and only if all of its Sylow subgroups are cyclic, so that for supersolvable groups all these conditions are equivalent. We also show that $G$ has property $\mathcal{D}$ whenever $N$ and $G/N$ do, provided $N$ is a normal Hall subgroup or the unique minimal normal subgroup of $G$. Applying this to Frobenius groups whose kernel is a minimal normal subgroup, we obtain infinitely many non-supersolvable groups with property $\mathcal{D}$, including the affine groups $\mathrm{AGL}(1,q)$ for every prime power $q$ that is not a prime. A computer search shows that there are exactly $88$ non-supersolvable groups with property $\mathcal{D}$ of order less than $200$, and that our results account for all of them.
\end{abstract}


\end{frontmatter}

\section{Introduction}\label{sec:intro}
It is well known that a finite group of order $n$ is cyclic if and only if it has a unique subgroup of order $d$ for each divisor $d$ of $n$. Put differently, the cyclic groups are the only finite groups in which distinct subgroups have distinct orders. In this paper we ask for the same thing, but only for \emph{normal} subgroups.

\begin{definition}\label{def:D}
A finite group $G$ is said to have \emph{property $\mathcal{D}$} if all normal subgroups of $G$ have distinct orders.
\end{definition}

Cyclic groups and simple groups clearly have property $\mathcal{D}$. An abelian group with property $\mathcal{D}$ is cyclic: every subgroup of an abelian group is normal, so distinct subgroups must have distinct orders. In fact, we will see in Lemma~\ref{lem:basic} that the same holds for nilpotent groups, so the interesting groups with property $\mathcal{D}$ are the non-nilpotent ones.

Since conjugate normal subgroups are equal, a group has property $\mathcal{D}$ if and only if its normal subgroups of equal order are conjugate. This puts $\mathcal{D}$ in a family of conditions that has been studied for a while: finite groups in which subgroups of the same order are conjugate (property~$\mathcal{B}$) \cite{vanderwall1}; in which abelian subgroups of equal order are conjugate (property~$\mathcal{A}$) \cite{sezer2}; in which cyclic subgroups of equal order are conjugate (property~$\mathcal{C}$) \cite{costantini-jabara}; in which supersolvable (respectively nilpotent) subgroups of equal order are conjugate (property~$\mathcal{SS}$, respectively $\mathcal{N}$) \cite{vanderwall2}; and in which cyclic $p$-subgroups of equal order are conjugate for every prime $p$ dividing the order of the group (property~$\mathcal{P}$) \cite{sezer1,vanderwall3}. These properties are related by
$$\mathcal{B}=\mathcal{SS}\Rightarrow \mathcal{N}\Rightarrow \mathcal{A}\Rightarrow \mathcal{C}\Rightarrow \mathcal{P}\Rightarrow \mathcal{D}.$$
Apart from the first equality and the last implication, this is clear from the definitions. The equality $\mathcal{B}=\mathcal{SS}$ is proved in \cite[Theorem 2.2]{vanderwall2}, and the last implication is Theorem~\ref{thm:P-implies-D} below. The implication cannot be reversed: $S_4$ has property $\mathcal{D}$ but not property $\mathcal{P}$ (Example~\ref{ex:S4}). So $\mathcal{D}$ is the weakest condition of this type studied so far, which is one more reason to look at these groups.

The paper is organized as follows. In Section~\ref{sec:basic} we collect some basic properties of groups with property $\mathcal{D}$. Section~\ref{sec:chain} proves $\mathcal{P}\Rightarrow\mathcal{D}$. In Section~\ref{sec:ss} we classify the supersolvable groups with property $\mathcal{D}$, showing that the whole chain above collapses in that case, and we give a criterion for solvable groups. Section~\ref{sec:ext} contains results that build new groups with property $\mathcal{D}$ from old ones, and Section~\ref{sec:frobenius} applies them to Frobenius groups and gives examples of non-supersolvable groups with property $\mathcal{D}$. Finally, Section~\ref{sec:conclusion} reports a computer search for small non-supersolvable groups with property $\mathcal{D}$ and lists some open problems.

\section{Basic properties}\label{sec:basic}
Throughout the paper all groups are finite, and $\mathbb{Z}_n$ denotes the cyclic group of order $n$. For a group $G$, we write $Z(G)$, $G'$, $\Phi(G)$ and $F(G)$ for its center, the derived subgroup, the Frattini subgroup and the Fitting subgroup, respectively. In this section we record the elementary consequences of property $\mathcal{D}$ that are used in the rest of the paper. The most useful of them is that the property is inherited by quotients; it is not inherited by subgroups, as the normal subgroup $\mathbb{Z}_2 \times \mathbb{Z}_2$ of $A_4$ shows.

\begin{lemma}\label{lem:basic}
Let $G$ be a finite group with property $\mathcal{D}$. Then:
\begin{enumerate}
\item Any quotient of $G$ also has property $\mathcal{D}$.
\item Every normal subgroup of $G$ is characteristic.
\item $Z(G)$ and $G/G'$ are cyclic.
\item If $G$ is nilpotent, then $G$ is cyclic.
\end{enumerate}
\end{lemma}
\begin{proof}
\begin{enumerate}
    \item Let $N\lhd G$. By the correspondence theorem, every normal subgroup of $G/N$ has the form $M/N$, where $M$ is a normal subgroup of $G$ containing $N$. Suppose $M_1/N$ and $M_2/N$ are normal subgroups of $G/N$ of the same order. Then $|M_1| = |M_2|$. As $G$ has property $\mathcal{D}$, $M_1 = M_2$, and hence $M_1/N = M_2/N$. So $G/N$ also has property $\mathcal{D}$.
    \item Let $N\lhd G$ and $\alpha\in \mathrm{Aut}(G)$. Then $\alpha(N)$ is a normal subgroup of $G$ with $|\alpha(N)|=|N|$. Hence $\alpha(N)=N$, i.e., $N$ is a characteristic subgroup.
    \item Every subgroup of $Z(G)$ is normal in $G$, so distinct subgroups of $Z(G)$ have distinct orders. As $Z(G)$ is abelian, it is cyclic. Since $G/G'$ is abelian, it is cyclic by (1).
    \item We first prove it for $p$-groups. Let $P$ be a $p$-group with property $\mathcal{D}$. By (1), $P/\Phi(P)$ has property $\mathcal{D}$, and by Burnside's basis theorem $P/\Phi(P)\cong \mathbb{Z}_p^{d}$ for some $d$. In particular $P/\Phi(P)$ is abelian with property $\mathcal{D}$, hence cyclic, so $d=1$. Thus $P$ is generated by a single element, i.e., $P$ is cyclic.

    Now let $G$ be nilpotent. Then it is the direct product of its Sylow subgroups $P_1, \ldots, P_k$. Each $P_i$ is isomorphic to a quotient of $G$ (divide by the product of the other Sylow subgroups), so $P_i$ has property $\mathcal{D}$ by (1) and is therefore cyclic. As a result, $G$ is cyclic.\qedhere
\end{enumerate}
\end{proof}

So a non-cyclic group with property $\mathcal{D}$ is never nilpotent, for example  $S_3$ and $A_4$. The next lemma is a standard fact about Sylow subgroups of normal subgroups. We state it separately because it will be used in the proofs of both Proposition \ref{prop:cyclic-sylow} and Theorem \ref{thm:P-implies-D}.

\begin{lemma}\label{lem:sylow}
Let $G$ be a finite group, $N\lhd G$, and let $P$ be a Sylow $p$-subgroup of $G$. Then $N\cap P$ is a Sylow $p$-subgroup of $N$. Consequently, if $P_i$ is a Sylow $p_i$-subgroup of $G$ for each prime $p_i$ dividing $|G|$, then $N$ is generated by the subgroups $N\cap P_i$.
\end{lemma}
\begin{proof}
Clearly $P\leq NP\leq G$. By the second isomorphism theorem, $NP/N \cong  P/(N \cap P)$, so $|NP| / |P| = |N| / |N \cap P|$. As $p\nmid [G:P]$ and $[G:P]=[G:NP] \cdot [NP:P]$, $p$ does not divide $[NP:P]=[N:N\cap P]$. Since $N\cap P\leq P$ is a $p$-group, $N \cap P$ is a Sylow $p$-subgroup of $N$.

For the second statement, the subgroup generated by the $N\cap P_i$ lies in $N$ and, by the first part, has order divisible by the full $p_i$-part of $|N|$ for every $i$. Hence it is all of $N$.
\end{proof}

\section{Property \texorpdfstring{$\mathcal{P}$}{P} implies property \texorpdfstring{$\mathcal{D}$}{D}}\label{sec:chain}
Recall that a group has property $\mathcal{P}$ if, for every prime $p$, any two of its cyclic $p$-subgroups of the same order are conjugate. In this section we prove the last implication in the chain from the introduction.

\begin{theorem}\label{thm:P-implies-D}
Let $G$ be a finite group with property $\mathcal{P}$. Then $G$ has property $\mathcal{D}$.
\end{theorem}
\begin{proof}
Let $N \lhd G$, let $P$ be a Sylow $p$-subgroup of $G$, and let $p^e$ be the largest order of an element of $N \cap P$. We claim that $$N \cap P = \{ x \in P : x^{p^e} = 1 \}.$$
The inclusion $\subseteq$ is clear. Conversely, let $x \in P$ have order $p^k$ with $k \leq e$. Since $N$ contains an element of order $p^e$, it contains a cyclic subgroup $Y$ of order $p^k$. As $G$ has property $\mathcal{P}$, there exists $g \in G$ with $\langle x \rangle = gYg^{-1}\leq gNg^{-1} = N$ (as $N\lhd G$). Hence $x \in N$, and the reverse inclusion holds. (This claim has also been proved in \cite[Lemma 2.1]{sezer1} in a different form. We give a separate proof to make the paper self-contained.)

Now let $N_1, N_2 \lhd G$ with $|N_1| = |N_2|$, and let $p^{e_i}$ be the largest order of an element of $N_i \cap P$. Without loss of generality, assume $e_1 \leq e_2$. Then, by the claim, $N_1 \cap P \subseteq N_2 \cap P$. By Lemma~\ref{lem:sylow}, $N_i \cap P$ is a Sylow $p$-subgroup of $N_i$, so $|N_1 \cap P| = |N_2 \cap P|$, and hence $N_1 \cap P = N_2 \cap P$. As this holds for a Sylow $p$-subgroup of $G$ for each prime $p$, and $N_i$ is generated by the subgroups $N_i \cap P$ (Lemma~\ref{lem:sylow}), we get $N_1 = N_2$. Thus $G$ has property $\mathcal{D}$.
\end{proof}

\begin{example}\label{ex:S4}
The symmetric group $S_4$ has property $\mathcal{D}$, since $1 < V_4 < A_4 < S_4$ are its only normal subgroups. It does not have property $\mathcal{P}$, since its subgroups of order $2$ are not all conjugate. So the implication $\mathcal{P}\Rightarrow\mathcal{D}$ is strict.
\end{example}

\section{Supersolvable and solvable groups}\label{sec:ss}
We now turn to the main result of the paper, the classification of supersolvable groups with property $\mathcal{D}$. One direction is easy and holds without any solvability assumption: a group all of whose Sylow subgroups are cyclic has property $\mathcal{D}$ (Proposition \ref{prop:cyclic-sylow}). The converse, Theorem \ref{thm:supersolvable}, is where supersolvability is used. As consequences, we show that all the conditions in the chain from the introduction coincide for supersolvable groups, and we describe when a solvable group with property $\mathcal{D}$ is supersolvable.

\begin{proposition}\label{prop:cyclic-sylow}
Let $G$ be a finite group in which every Sylow subgroup is cyclic. Then $G$ has property $\mathcal{D}$.
\end{proposition}
\begin{proof}
Let $N_1, N_2\lhd G$ with $|N_1|=|N_2|$. By Lemma~\ref{lem:sylow}, for each Sylow $p$-subgroup $P$ of $G$, the subgroups $N_1\cap P$ and $N_2\cap P$ have the same order, and both lie in the cyclic group $P$. So $N_1\cap P=N_2\cap P$ for every Sylow subgroup $P$ of $G$. Now $N_1$ and $N_2$ are generated by the subgroups $N_1\cap P$ and $N_2\cap P$ respectively (Lemma~\ref{lem:sylow}), so $N_1=N_2$. Thus $G$ has property $\mathcal{D}$.
\end{proof}

For the converse we need two standard facts about supersolvable groups.

\begin{proposition}\label{prop:ss-facts}
Let $G$ be a finite supersolvable group. Then:
\begin{enumerate}[label=(\arabic*)]
    \item If $p$ is the largest prime dividing $|G|$, then the Sylow $p$-subgroup of $G$ is normal in $G$. \cite[Theorem 9.1 (c)]{huppert1}
    \item If $N$ is a non-trivial normal subgroup of $G$, then $N$ contains a normal subgroup of $G$ of prime order. (This follows from the fact that every minimal normal subgroup of a finite supersolvable group has prime order.)
\end{enumerate}
\end{proposition}

\begin{theorem}\label{thm:supersolvable}
Let $G$ be a finite supersolvable group. Then $G$ has property $\mathcal{D}$ if and only if every Sylow subgroup of $G$ is cyclic.
\end{theorem}
\begin{proof}
If every Sylow subgroup of $G$ is cyclic, then $G$ has property $\mathcal{D}$ by Proposition \ref{prop:cyclic-sylow}.

Conversely, let $G$ be supersolvable with property $\mathcal{D}$. We prove that every Sylow subgroup is cyclic by strong induction on $|G|$. Let $p$ be the largest prime dividing $|G|$ and $P$ the Sylow $p$-subgroup of $G$. Then by Proposition~\ref{prop:ss-facts}(1), $P$ is normal in $G$. Thus $G/P$ is supersolvable (quotients of supersolvable groups are supersolvable) with property $\mathcal{D}$ (Lemma~\ref{lem:basic}(1)), and has strictly smaller order than $|G|$. So by the induction hypothesis its Sylow subgroups are cyclic. If $R$ is a Sylow $r$-subgroup of $G$ with $r\neq p$, then $R \cap P$ is trivial and $R \cong RP/P$, which is a Sylow subgroup of $G/P$. So $R$ is cyclic. It remains to show that $P$ is cyclic.

Suppose, to the contrary, that $P$ is not cyclic. As $\Phi(P)$ is characteristic in $P$ and $P \lhd G$, we have $\Phi(P) \lhd G$. Put $\widehat{G} = G/\Phi(P)$ and $V = P/\Phi(P)$. Then $\widehat{G}$ is supersolvable and has property $\mathcal{D}$, and $V \lhd \widehat{G}$. By the Burnside basis theorem, $V$ is an elementary abelian $p$-group, and $|V| \geq p^2$ since $P$ is not cyclic. Moreover $[\widehat{G}:V] = [G:P]$ is not divisible by $p$.

By Proposition~\ref{prop:ss-facts}(2), $V$ contains a subgroup $U \lhd \widehat{G}$ with $|U| = p$. As $V$ is an elementary abelian $p$-group, $U$ has a complement in $V$. Since $U$ is an abelian normal subgroup of $\widehat{G}$, $U \leq V \leq \widehat{G}$, and $[\widehat{G}:V]$ is coprime to $|U| = p$, it follows from Gasch\"utz's theorem \cite[Main Theorem 17.4 (a)]{huppert1} that $U$ has a complement $H$ in $\widehat{G}$, i.e., $\widehat{G} = UH$ and $U \cap H = 1$.

Put $V_1 = V \cap H$. As $U \leq V$, we have $V = V \cap UH = U(V \cap H) = UV_1$, and $U \cap V_1 = 1$. Hence $|V_1| = |V|/p \geq p$, so $V_1$ is non-trivial.

\smallskip
\noindent\textit{Claim:} $V_1 \lhd \widehat{G}$.\\
\noindent\textit{Proof of the claim.} Since $\widehat{G} = UH$, every element of $\widehat{G}$ is a product $uh$ with $u \in U$ and $h \in H$. So it is enough to show that conjugation by elements of $H$ and by elements of $U$ both map $V_1$ to itself.

\emph{Conjugation by elements of $U$:} Let $u \in U$ and $x \in V_1$. Both $u,x \in V$ and $V$ is abelian, so $ux=xu$, i.e., $uxu^{-1}=x\in V_1$.

\emph{Conjugation by elements of $H$:} Let $h \in H$ and $x \in V_1 = V \cap H$. As $x \in V$ and $V \lhd \widehat{G}$, we have $hxh^{-1} \in V$. On the other hand, $x \in H$ and $h \in H$, so $hxh^{-1} \in H$. Thus $hxh^{-1}\in V_1$, and the claim follows.

\smallskip
Again by Proposition \ref{prop:ss-facts}(2), $V_1$ contains a subgroup $U' \lhd \widehat{G}$ with $|U'| = p$. Since $U \cap V_1$ is trivial, $U' \neq U$. Thus $U$ and $U'$ are distinct normal subgroups of $\widehat{G}$ of the same order $p$, which contradicts the fact that $\widehat{G}$ has property $\mathcal{D}$. Hence $P$ is cyclic, and this completes the proof.
\end{proof}

Combined with the chain from the introduction, Theorem~\ref{thm:supersolvable} shows that for supersolvable groups all the conditions in the family collapse into one.

\begin{corollary}\label{cor:equivalent}
Let $G$ be a finite supersolvable group. Then the following are equivalent: $G$ has property $\mathcal{B}$; $G$ has property $\mathcal{N}$; $G$ has property $\mathcal{A}$; $G$ has property $\mathcal{C}$; $G$ has property $\mathcal{P}$; $G$ has property $\mathcal{D}$; every Sylow subgroup of $G$ is cyclic.
\end{corollary}
\begin{proof}
By the chain $\mathcal{B} \Rightarrow \mathcal{N} \Rightarrow \mathcal{A} \Rightarrow \mathcal{C} \Rightarrow \mathcal{P} \Rightarrow \mathcal{D}$ (the last step being Theorem \ref{thm:P-implies-D}), it suffices to show that property $\mathcal{D}$ implies property $\mathcal{B}$. If $G$ has property $\mathcal{D}$, then all Sylow subgroups of $G$ are cyclic by Theorem \ref{thm:supersolvable}, and such a group has property $\mathcal{B}$ by \cite[p.~450]{vanderwall1}.
\end{proof}

For solvable groups with property $\mathcal{D}$, supersolvability can be read off from the Fitting subgroup.

\begin{corollary}\label{cor:fitting}
Let $G$ be a finite solvable group with property $\mathcal{D}$. Then the following are equivalent:
\begin{enumerate}
\item $G$ is supersolvable;
\item every Sylow subgroup of $G$ is cyclic;
\item $F(G)$ is cyclic.
\end{enumerate}
\end{corollary}
\begin{proof}
(1)$\Rightarrow$(2) is Theorem \ref{thm:supersolvable}. For (2)$\Rightarrow$(3), if every Sylow subgroup of $G$ is cyclic, then each $O_p(G)$ is cyclic (as it is contained in some Sylow $p$-subgroup), and so is their direct product $F(G)$.

For (3)$\Rightarrow$(1), let $F = F(G)$ be cyclic. As $G$ is solvable, $C_G(F) \leq F$. Also, as $F$ is abelian, $F\leq C_G(F)$. Thus $C_G(F)=F$. So, by the N/C-theorem, $G/F$ embeds in $\mathrm{Aut}(F)$, which is abelian. Hence $G' \leq F$ is cyclic, and $G/G'$ is cyclic by Lemma \ref{lem:basic}(3). Thus $G$ is metacyclic, hence supersolvable.
\end{proof}

\begin{remark}\label{rem:fitting}
It follows from Corollary \ref{cor:fitting} that a solvable group with property $\mathcal{D}$ which is not supersolvable has a non-cyclic Fitting subgroup.
\end{remark}

\section{Extensions and direct products}\label{sec:ext}
By Theorem \ref{thm:supersolvable}, a non-supersolvable group with property $\mathcal{D}$ must be built in a different way, and in this section we give tools for doing so. By Lemma \ref{lem:basic}(1), if $G$ has property $\mathcal{D}$ then so does $G/N$ for every normal subgroup $N$. The converse question is whether property $\mathcal{D}$ for $N$ and $G/N$ forces it for $G$. In general it does not: the dihedral group of order $12$ has a cyclic normal subgroup $N$ of order $6$ with $G/N \cong \mathbb{Z}_2$, but it has three normal subgroups of order $6$. We show that the answer is positive when $N$ is a normal Hall subgroup or the unique minimal normal subgroup of $G$, and then relax the second condition.

\begin{theorem}\label{thm:extension}
Let $G$ be a finite group such that $G/N$ and $N$ both have property $\mathcal{D}$, where $N$ is either a normal Hall subgroup of $G$ or the unique minimal normal subgroup of $G$. Then $G$ also has property $\mathcal{D}$.
\end{theorem}
\begin{proof}
\textit{Case 1: $N$ is a normal Hall subgroup of $G$.} Let $K_1, K_2 \lhd G$ with $|K_1| = |K_2|$. For $i = 1,2$ we have
$|K_i| = |K_i \cap N| \cdot |K_iN/N|$, where the first factor divides $|N|$ and the second divides $[G:N]$. As these are coprime, $|K_1 \cap N| = |K_2 \cap N|$ and $|K_1N/N| = |K_2N/N|$. Since $K_i \cap N \lhd N$ and $K_iN/N \lhd G/N$, property $\mathcal{D}$ of $N$ and of $G/N$ gives $K_1 \cap N = K_2 \cap N = L$ (say) and $K_1N = K_2N = M$ (say).

Let $x \in K_1$ and put $k = [M:N]$. Since
$K_1/L = K_1/(K_1 \cap N) \cong K_1N/N = M/N$ has order $k$, we have $x^k \in L \leq K_2$. On the other hand, $x \in M = K_2N$ and
$M/K_2 = K_2N/K_2 \cong N/(K_2 \cap N) = N/L$, so $x^{[N:L]} \in K_2$. As $k$ divides $[G:N]$ and $[N:L]$ divides $|N|$, the integers $k$ and $[N:L]$ are coprime, so there are integers $a, b$ with $ak + b[N:L] = 1$. Hence $$x = (x^k)^a (x^{[N:L]})^b \in K_2 .$$
Thus $K_1 \leq K_2$. Similarly $K_2\leq K_1$, and hence $K_1=K_2$. So $G$ has property $\mathcal{D}$.

\smallskip
\textit{Case 2: $N$ is the unique minimal normal subgroup of $G$.} Let $K_1, K_2 \lhd G$ with $|K_1| = |K_2|$. If $K_1 = 1$, then $K_2 = 1$ and we are done. So assume that $K_1$ and $K_2$ are non-trivial. As every non-trivial normal subgroup contains a minimal normal subgroup, we have $N \leq K_1,K_2$. Thus $K_1/N$ and $K_2/N$ are normal subgroups of $G/N$ of the same order. As $G/N$ has property $\mathcal{D}$, we get $K_1/N = K_2/N$, and therefore $K_1 = K_2$. So $G$ has property $\mathcal{D}$.
\end{proof}

\begin{remark}\label{rem:extension}
In the first case of Theorem~\ref{thm:extension}, we can replace ``$N$ has property $\mathcal{D}$'' by the slightly weaker condition ``the normal subgroups of $G$ contained in $N$ have distinct orders'', without changing the proof. In the second case, i.e., when $N$ is the unique minimal normal subgroup of $G$, we do not need $N$ to have property $\mathcal{D}$ at all.
\end{remark}

The first case of Theorem \ref{thm:extension} gives the following description of direct products of groups of coprime orders.

\begin{corollary}\label{cor:coprime}
Let $G$ and $H$ be finite groups with $\gcd(|G|,|H|) = 1$. Then $G \times H$ has property $\mathcal{D}$ if and only if both $G$ and $H$ have property $\mathcal{D}$.
\end{corollary}
\begin{proof}
If $G \times H$ has property $\mathcal{D}$, then by Lemma~\ref{lem:basic}(1), both $G$ and $H$ have property $\mathcal{D}$ (they are quotients of $G\times H$).

Conversely, let $G$ and $H$ have property $\mathcal{D}$. As $H$ is a normal Hall subgroup of $G\times H$ with quotient $(G\times H)/H\cong G$, Case 1 of Theorem~\ref{thm:extension} shows that $G \times H$ has property $\mathcal{D}$.
\end{proof}

The second case of Theorem \ref{thm:extension} can be pushed a little further, by allowing other minimal normal subgroups as long as they are controlled.

\begin{theorem}\label{thm:minimal-normal}
Let $G$ be a finite group and $N$ be a minimal normal subgroup of $G$ such that $G/N$ has property $\mathcal{D}$. If the normal subgroups $K$ of $G$ with $K \cap N = 1$ have pairwise distinct orders, none divisible by $|N|$, then $G$ has property $\mathcal{D}$.
\end{theorem}
\begin{proof}
Let $K_1, K_2 \lhd G$ with $|K_1| = |K_2|$. For $i=1,2$, the subgroup $K_i \cap N$ is normal in $G$ and contained in $N$, so by minimality of $N$ either $N \leq K_i$ or $K_i \cap N = 1$.

Suppose that $N \leq K_1$ but $K_2 \cap N = 1$. Then $|N|$ divides $|K_1| = |K_2|$, which contradicts the assumption that the order of a normal subgroup meeting $N$ trivially is not divisible by $|N|$. By symmetry, either both $K_1,K_2$ contain $N$, or both meet $N$ trivially.

If both contain $N$, then $K_1/N$ and $K_2/N$ are normal subgroups of $G/N$ of the same order. As $G/N$ has property $\mathcal{D}$, $K_1/N=K_2/N$, and so $K_1=K_2$. If both meet $N$ trivially, then $K_1=K_2$ by the assumption that such subgroups have pairwise distinct orders. Hence $G$ has property $\mathcal{D}$.
\end{proof}

\begin{remark}\label{rem:minimal-normal}
Theorem \ref{thm:minimal-normal} applies to groups that are out of reach of Theorem \ref{thm:extension}, even in the weaker form of Remark \ref{rem:extension}. 

 For example, let $G = Z \times A$ with $Z \cong \mathbb{Z}_2$ and $A \cong A_5$. The proper non-trivial normal subgroups of $G$ are $Z$ and $A$. Both are minimal normal subgroups, so $G$ has no unique minimal normal subgroup, and neither is a Hall subgroup, as $[G:Z] = 60$ and $[G:A] = 2$. Hence Theorem \ref{thm:extension} does not apply. Now take $N = A$ in Theorem \ref{thm:minimal-normal}. Then $G/N \cong \mathbb{Z}_2$ has property $\mathcal{D}$. If $K \lhd G$ and $K \cap A = 1$, then for $k \in K$ and $a \in A$ the element $kak^{-1}a^{-1}$ lies in $K \cap A = 1$, as both $K$ and $A$ are normal in $G$. So $K$ commutes elementwise with $A$, i.e., $K \leq C_G(A) = Z$. Thus the normal subgroups of $G$ meeting $A$ trivially are $1$ and $Z$, of orders $1$ and $2$, neither divisible by $|A| = 60$. Hence $G$ has property $\mathcal{D}$.
\end{remark}

\section{Frobenius groups and examples}\label{sec:frobenius}
We close by applying the results of Section \ref{sec:ext} to Frobenius groups. In a Frobenius group every normal subgroup either lies in the Frobenius kernel or contains it, which is exactly the situation in which property $\mathcal{D}$ can be read off from the kernel and the complement. This gives an infinite supply of groups with property $\mathcal{D}$ that are not supersolvable, and so are not covered by Theorem \ref{thm:supersolvable}.

\begin{theorem}\label{thm:frobenius}
Let $G$ be a Frobenius group with Frobenius kernel $F$ and Frobenius complement $H$. Suppose that $F$ is a minimal normal subgroup of $G$. Then $G$ has property $\mathcal{D}$ if and only if $H$ has property $\mathcal{D}$. In particular, $G$ has property $\mathcal{D}$ whenever every Sylow subgroup of $H$ is cyclic.
\end{theorem}
\begin{proof}
We have $G = FH$ and $F \cap H = 1$, so $G/F \cong H$. If $G$ has property $\mathcal{D}$, then so does $H$ by Lemma \ref{lem:basic}(1).

Conversely, suppose that $H$ has property $\mathcal{D}$. Let $N$ be a non-trivial normal subgroup of $G$. By \cite[Chapter V, Theorem 8.16]{huppert1}, either $N \leq F$ or $F < N$. In the first case $N = F$, as $F$ is a minimal normal subgroup of $G$. Thus every non-trivial normal subgroup of $G$ contains $F$.

Now let $K_1, K_2 \lhd G$ with $|K_1| = |K_2|$. If $K_1 = 1$, then $K_2 = 1$. Otherwise $K_1$ and $K_2$ both contain $F$, so $K_1/F$ and $K_2/F$ are normal subgroups of $G/F \cong H$ of the same order. As $H$ has property $\mathcal{D}$, we get $K_1/F = K_2/F$, and hence $K_1 = K_2$. Thus $G$ has property $\mathcal{D}$.

The last statement follows from Proposition \ref{prop:cyclic-sylow}, by which a group with all Sylow subgroups cyclic has property $\mathcal{D}$.
\end{proof}

\begin{example}\label{ex:frobenius}
The following Frobenius groups have a minimal normal Frobenius kernel and a cyclic Frobenius complement, and hence have property $\mathcal{D}$ by Theorem \ref{thm:frobenius}:
$$A_4 \cong \mathbb{Z}_2^2 \rtimes \mathbb{Z}_3, \quad
\mathbb{Z}_3^2 \rtimes \mathbb{Z}_4, \quad
\mathbb{Z}_2^3 \rtimes \mathbb{Z}_7, \quad
\mathbb{Z}_3^2 \rtimes \mathbb{Z}_8, \quad
\mathbb{Z}_5^2 \rtimes \mathbb{Z}_3, \quad
\mathbb{Z}_2^4 \rtimes \mathbb{Z}_5 .$$
None of these groups is supersolvable, since the kernel is a minimal normal subgroup of non-prime order. More generally, for a prime power $q$, the group $\mathrm{AGL}(1,q) \cong \mathbb{F}_q \rtimes \mathbb{F}_q^*$ has property $\mathcal{D}$; it is not supersolvable as soon as $q$ is not a prime.
\end{example}

\section{Concluding remarks and open problems}\label{sec:conclusion}

We introduced property $\mathcal{D}$ as the weakest member of a family of conjugacy
conditions on subgroups of equal order, and determined the supersolvable groups with this property: they are exactly the groups all of whose Sylow subgroups are cyclic (Theorem \ref{thm:supersolvable}). Outside the supersolvable case the picture is much less rigid. Theorems \ref{thm:extension}, \ref{thm:minimal-normal} and \ref{thm:frobenius} produce many non-supersolvable examples, but they do not amount to a classification. We end with some observations and questions that we believe deserve further study.

We first note that property $\mathcal{D}$ for \emph{all} groups of a given order is easy to describe: every group of order $n$ has property $\mathcal{D}$ if and only if $n$ is squarefree. Indeed, if $n$ is squarefree, then every Sylow subgroup of a group of order $n$ has prime order, and Proposition \ref{prop:cyclic-sylow} applies; if $p^2$ divides $n$, then $\mathbb{Z}_p \times \mathbb{Z}_p \times \mathbb{Z}_{n/p^2}$ is a non-cyclic abelian group of order $n$. The existence question is more delicate.

\begin{problem}
Determine the positive integers $n$ for which there exists a non-cyclic group of order $n$ with property $\mathcal{D}$. By Lemma \ref{lem:basic}(4), a necessary condition is that there exists a non-nilpotent group of order $n$. Is this condition sufficient?
\end{problem}

We used \textsf{SageMath} \cite{sagemath}, which accesses the SmallGroups library through \textsf{GAP} \cite{gap4}, to list all groups with property $\mathcal{D}$ of order less than $200$. Exactly $30$ non-supersolvable groups of order less than $100$, and $88$ of order less than $200$, have property $\mathcal{D}$; all of them are solvable except $A_5$, $\mathrm{SL}(2,5)$, $S_5$, $\mathbb{Z}_2 \times A_5$, $\mathrm{PSL}(3,2)$ and $\mathrm{GL}(2,4)$. 

We also verified by computer that each of these $88$ groups satisfies the hypotheses of at least one of Proposition \ref{prop:cyclic-sylow} and Theorems \ref{thm:extension} (in the form of Remark \ref{rem:extension}), \ref{thm:minimal-normal} and \ref{thm:frobenius}, where the quotients and subgroups appearing in those hypotheses are handled by the same results. So these results account for every non-supersolvable group with property $\mathcal{D}$ in this range, which motivates the following.

\begin{problem}
Classify the finite solvable groups with property $\mathcal{D}$. By Corollary
\ref{cor:fitting}, it remains to treat those with non-cyclic Fitting subgroup.
\end{problem}

Finally, the non-solvable case is wide open. Every non-abelian simple group has property $\mathcal{D}$, and if $T$ is a non-abelian simple group and $T \leq G \leq
\mathrm{Aut}(T)$, then $T$ is the unique minimal normal subgroup of $G$, so by Theorem \ref{thm:extension} and Remark \ref{rem:extension}, $G$ has property $\mathcal{D}$ if and only if $G/T$ does. Beyond such cases little is known.

\begin{problem}
Describe the non-solvable groups with property $\mathcal{D}$. In particular, determine when a direct product of non-abelian simple groups, and when a perfect group, has property $\mathcal{D}$.
\end{problem}

\section*{Statements and Declarations}
The author acknowledges the funding of DST-FIST Sanction no. SR/FST/MS-I/2019/41, Govt. of India.

\subsection*{Declaration of use of generative AI}
The author used Claude (Anthropic) (Free version) to cross-check the computations reported in Section \ref{sec:conclusion}, to identify related literature, and to improve the language and organization of the manuscript.

\subsection*{Data Availability Statement}
No datasets were generated or analyzed during the current study. The \textsf{SageMath} code used for the computations reported in Section \ref{sec:conclusion} is available from the author on request.

\subsection*{Competing Interests}
The author has no competing interests to declare that are relevant to the content of this article.

\end{document}